\documentclass[a4paper,11pt]{article}
\usepackage[utf8]{inputenc} 
\usepackage{amsmath,amsfonts,amssymb,graphics,epsfig,color}

\usepackage[left=2cm,top=3cm,right=2cm,bottom=3cm,bindingoffset=0.5cm]{geometry}
\usepackage[english]{babel}
\usepackage{verbatim}
\usepackage{amscd,enumerate}
\usepackage{epstopdf}
\usepackage{graphicx}
\usepackage{multicol}
\usepackage{appendix}
\usepackage{mathrsfs}
\usepackage{amsmath}
\usepackage{amsthm}
\usepackage{cite}
\newcommand\be{\begin{equation}}
\newcommand\ee{\end{equation}}

\title{Hybrid Ermakov-Ray-Reid/Painlevé II Symmetry Reduction: Application to a Class of Moving Boundary Problems.}
\author{
Colin Rogers $^{1}$, Adriana C. Briozzo $^{2}$\\
\small {{$^1$} School of Mathematics and Statistics, The University of New South Wales},\\
\small {Sydney NSW 2052, Australia}\\
\small {{$^2$} Depto. Matem\'atica, FCE, Univ. Austral, Paraguay 1950-CONICET} \\
\small {S2000FZF Rosario, Argentina.}\\
\small{Email: c.rogers@unsw.edu.au,  abriozzo@austral.edu.ar}
}

\begin{document}

\maketitle

\begin{abstract}
Here, a class of nonlinear moving boundary problems for a novel extension of a two-component mKdV system is shown to admit exact solution via application of a hybrid Ermakov-Ray-Reid /
Painlevé II symmetry ansatz.The mKdV system has its genesis in a reduction of a coupled nonlinear NLS system incorporating 
deBroglie - Bohm potential terms.

\end{abstract}

\section{\textbf{Introduction}} 
Ermakov-type systems with genesis in the classical work of \cite{ermakov1880} have diverse physical applications in both physics and nonlinear continuum mechanics \cite{rogers2018b}. The original Ermakov equation together with its admitted nonlinear superposition principle has application notably in the analysis of boundary value problems associated with the large amplitude oscillation of thin-walled cylinders composed of hyperelastic Mooney-Rivlin material subject to a range of boundary loadings \cite{shahinpoor1971, rogers1989b}. In addition, the classical Ermakov equation occurs in the analysis of physical pulsrodon phenomena in oceanographic warm-core eddy evolution theory \cite{rogers1989c}. The canonical nonlinear superposition principle admitted by this Ermakov equation can be derived via Lie group invariance considerations as in the hydrodynamic context of \cite{rogers1989a}.
Nonlinear two-component coupled systems of Ermakov-Ray-Reid type were originally introduced in \cite{ray1980,reid1980} with extension to 2+1-dimensional Ermakov systems in \cite{rogers1993}. Ermakov-Ray-Reid systems have diverse physical applications such as in nonlinear optics \cite{rogers2010a,rogers2012}, magnetogasdynamics \cite{rogers2010b}, rotating shallow water hydrodynamics \cite{rogers2010c}, elasticity \cite{rogers2013}, 2+1- dimensional Madelung hydrodynamic system analysis \cite{rogers2011} and many-body theory \cite{rogers2018}.

In \cite{rogers1996}, multicomponent Ermakov-type systems were derived via symmetry reduction of a multilayer hydrodynamic system. Application was subsequently made to the analysis of non autonomous many-body theory \cite{rogers2018}.

Hybrid Ermakov-Painlevé II systems were originally derived in \cite{rogers2014} in the context of wave packet representations admitted by multi-dimensional coupled nonlinear Schrodinger systems incorporations  de Broglie-Bohm potential. Therein, the canonical base Ermakov-Painlevé II euation was derived in the analysis of transverse wave propagation in a generalised Mooney-Rivlin hyperelastic material. It has later been derived in a range of physical applications, such as cold plasma physics \cite{rogers2018a}, Korteweg capillarity theory \cite{rogers2017c} and in the analysis of Dirichlet boundary value problems for the classical Nernst-Planck electrolysis system \cite{amster2015}.

Integrable structure underlying Ermakov-Painlevé II systems has been detailed in \cite{rogers2016b}. Therein, an admitted invariant was applied in an algorithmic solution procedure involving a Painlevé II connection. Iterative application of a Bäcklund transformation due to Lukashevich \cite{lukashovich1971} 
for Painlevé II was used therein. An important link between the canonical single component Ermakov-Painlevé II equation and the integrable Painlevé XXXIV equation was established.

Moving boundary problems of Stefan-type which admit exact solution via Painlevé II symmetry reduction was conducted for the canonical solitonic Dym equation (qv\cite{vassiliou2001}) and its reciprocal associates in \cite{rogers2015}. The physical motivation resided in the classical Saffman-Taylor model. The latter describes the percolation of a liquid into a porous medium or Hele-Shaw cell \cite{fokas1998}. In \cite{schief1999}, an integrable extension of the Dym equation was shown to be embedded in a novel geometric solitonic system descriptive of torsion evolution asssociated with the spatial binormal motion of inextensible curves. A reciprocal transformation \cite{kingston1982} was applied to established equivalence of this torsion system to the m$^{2}$KdV solitonic equation originally derived in \cite{fokas1980} and later, independently in \cite{calogero1985}. In terms of physical applications, the extended Dym equation occurs, importantly in the description of unidirectional dispersive shallow water propagation with novel associated peaked solitonic phenomena \cite{camassa1993}. In \cite{rogers2017b} certain moving boundary problems of Stefan-type for the extended Dym equation were shown to admit exact solution via Painlevé symmetry reduction. Iterated action of the Bäcklund transformation for Painlevé II was applied to extend the range of such solvable moving boundary problems. In addition, applications of reciprocal-type transformations to determine associated solvable moving boundary problems was detailed.

Moving boundary problems in 1+1- dimensional soliton theory which are exactly solvable via Painlevé II symmetry reduction have been analysed in \cite{rogers2016,rogers2022,rogers2023,rogers2025a,rogers2025b}. 
In \cite{rogers2025d}, a novel extension of the solitonic mKdV equation has recently been shown to be amenable to Ermakov-Painlevé II symetry reduction and application thereby made to solve a class of associated moving boundary problem.

Here, an extended mKdV coupled two-component system is shown to admit novel symmetry reduction to a hybrid Ermakov-Ray-Reid/Painlevé II system as originally set down in \cite{rogers2016b}. A class of moving boundary problems for the coupled mKdV-type system is shown thereby is admit exact solution. This class involves a region bounded by a pair of moving curves. Moving boundary problems for two-component coupled systems with this characteristic occur in modern soliton theory in connection with the resonant nonlinear Schrödinger equation \cite{flavin2008}. The latter and the resonant Davey-Stewartaon system have physical applications notably in capillarity theory \cite{rogers1999,rogers2009} and cold plasma physics \cite{lee2007}.

\section{\textbf{A Ermakov-Ray-Reid /Painlevé II System: Symmetry Reduction}}
Here, a novel extended mKdV-type two-component coupled system is introduced, namely

\begin{equation} \label{2.1}
\begin{split}
u_{t}+u_{xxx}+\alpha[(u^{2}+v^{2})u]_x+\lambda (t+a)^{\mu}\left[\frac{1}{u^{2}v}S(v/u)\right]_x &=0 \\[2ex]
v_{t}+v_{xxx}+\alpha[(u^{2}+v^{2})v]_x+\lambda (t+a)^{\mu}\left[\frac{1}{v^{2}u}T(u/v)\right]_x &=0, \quad \lambda, \mu \in \mathbb{R}
\end{split}
\end{equation}

The latter is shown to admit symmetry reduction to a hybrid Ermakov-Ray-Reid /Painlevé II system. 

\begin{center}
    Symmetry Reduction Ansätze
\end{center}

Representations 
\begin{equation}\label{2.2}
    u=(t+a)^{m} \Phi\bigg(\frac{x}{(t+a)^{n}}\bigg)\quad\quad
\end{equation}

\begin{equation}\label{2.3}
    v=(t+a)^{m} \Psi\bigg(\frac{x}{(t+a)^{n}}\bigg)\quad\quad
\end{equation}
are now inserted into the system \eqref{2.1}. In the case of $\eqref{2.1}_1$ this yields 

\begin{equation}\label{2.4}
    m(t+a)^{m-1}\Phi-n\xi(t+a)^{m-1}\Phi'+(t+a)^{m-3n}\Phi'''+\alpha[(\Phi^{2}+\Psi^{2})\Phi]'(t+a)^{3m-n}
\end{equation} 

\[+\lambda (t+a)^{\mu-3m-n}\left[\frac{1}{\Phi^{2}\Psi}S(\Psi/\Phi)\right]'=0,\quad (\xi=x/(t+a)^{n})\]
whence

\begin{equation}\label{2.5}
   \Psi'''+(m+n)\Phi-n(\xi\Phi)'+\alpha[(\Phi^{2}+\Psi^{2})\Phi]'+\lambda\left[\frac{1}{\Phi^{2}\Psi}S(\Psi/\Phi)\right]'=0
\end{equation}
together with $m=-\frac{1}{3}$, $n=\frac{1}{3}$and $\mu=-2.$ 
Thus
\begin{equation}\label{2.6}
   \Phi'''-(1/3)(\xi\Phi)'+\alpha[(\Phi^{2}+\Psi^{2})\Phi]'+\lambda\left[\frac{1}{\Phi^{2}\Psi}S(\Psi/\Phi)\right]'=0
\end{equation}and in a similar manner, insertion of the symmetry representation \eqref{2.3} in $\eqref{2.1}_2$ yields
\begin{equation}\label{2.7}
   \Psi'''-(1/3)(\xi\Psi)'+\alpha[(\Phi^{2}+\Psi^{2})\Psi]'+\lambda\left[\frac{1}{\Psi^{2}\Phi}T(\Phi/\Psi)\right]'=0
\end{equation}

On integration in turn of \eqref{2.6} and \eqref{2.7} there results 

\begin{equation}\label{2.8}
   \Phi''-(1/3)\xi\Phi+\alpha(\Phi^{2}+\Psi^{2})\Phi+\lambda\left[\frac{1}{\Phi^{2}\Psi}S(\Psi/\Phi)\right]=\alpha_I,
\end{equation}
\begin{equation}\label{2.9}
   \Psi''-(1/3)\xi\Psi+\alpha[(\Phi^{2}+\Psi^{2})\Psi]+\lambda\left[\frac{1}{\Psi^{2}\Phi}T(\Phi/\Psi)\right]=\alpha_{II}
\end{equation}
This with $\alpha_I=\alpha_{II}=0$ constitutes a hybrid Ermakov-Ray-Reid/Painlevé II system of the type originally introduced in \cite{rogers2016b}.

\section{A class of Moving Boundary Problems}
Here, moving boundary problems of Stefan-type are detailed for which the nonlinear coupled system \eqref{2.1} is amenable to exact solution via Ermakov-Ray-Reid/Painlevé II symmetry reduction.

The moving boundary problems under consideration concern a region bounded by $x=\Sigma_1(t)=\gamma_1(t+a)^{1/3}$ and $x=\Sigma_2(t)=\gamma_2(t+a)^{1/3}$ with $\gamma_1<\gamma_2$, $t>0$.
\begin{center}
    Boundary Conditions
\end{center}
\begin{equation} \label{3.1}
u_{xx}+\alpha[(u^{2}+v^{2})u]+\lambda (t+a)^{-2}\left[\frac{1}{u^{2}v}S(v/u)\right]=L_m
\Sigma_1^{i}\dot{\Sigma}_1, \quad \textbf{on}\quad x=\Sigma_1(t), t>0
\end{equation}

\begin{equation}\label{3.2}
v_{xx}+\alpha[(u^{2}+v^{2})v]+\lambda (t+a)^{-2}\left[\frac{1}{v^{2}u}T(u/v)\right]=M_m \Sigma_2^{i}\dot{\Sigma}_2, \quad \textbf{on}\quad x=\Sigma_2(t), t>0,
\end{equation}
together with
\begin{equation}\label{3.3}
    u=P_m \Sigma_1^{j}\quad \textbf{on}\quad x=\Sigma_1(t), t>0,
\end{equation}
\begin{equation}\label{3.4}
    v=R_m \Sigma_2^{l}\quad \textbf{on}\quad x=\Sigma_2 (t), t>0.
\end{equation}
The preceding constraints are now addressed in turn: 

I 
\[
u_{xx}+\alpha[(u^{2}+v^{2})u]+\lambda (t+a)^{-2}\left[\frac{1}{u^{2}v}S(v/u)\right]=L_m
\Sigma_1^{i}\dot{\Sigma}_1, \quad \textbf{on}\quad x=\Sigma_1(t), t>0
\]
Insertion of the symmetry ansatz \eqref{2.2} into the preceding yields 

    \begin{equation}\label{3.5}
\Phi''+\alpha(\Phi^{2}+\Psi^{2})\Phi+\lambda\left[\frac{1}{\Phi^{2}\Psi}S(\Psi/\Phi)\right]|_{\xi=\gamma_1}=(t+a)^{2n-m}L_m\Sigma_1^{i}\dot{\Sigma}_1
\end{equation}
whence $i=-1$ together with
 \begin{equation}\label{3.6}
\Phi''+\alpha(\Phi^{2}+\Psi^{2})\Phi+\lambda\left[\frac{1}{\Phi^{2}\Psi}S(\Psi/\Phi)\right]|_{\xi=\gamma_1}=1/3 L_m
\end{equation}
But \eqref{2.8} with $ \alpha_I=\alpha_{II}=0$ at $\xi=\gamma_1$ yields 
\begin{equation}\label{3.7}
\Phi''-(1/3) \xi \Phi+\alpha(\Phi^{2}+\Psi^{2})\Phi+\lambda\left[\frac{1}{\Phi^{2}\Psi}S(\Psi/\Phi)\right]|_{\xi=\gamma_1}=0
\end{equation}
which combined with \eqref{3.6} shows that
\begin{equation}\label{3.8}
    L_m=\gamma_1 \Phi(\gamma_1)
\end{equation}

II
\[
v_{xx}+\alpha[(u^{2}+v^{2})v]+\lambda (t+a)^{-2}\left[\frac{1}{v^{2}u}T(u/v)\right]=M_m
\Sigma_2^{i}\dot{\Sigma}_2, \quad \textbf{on}\quad x=\Sigma_2(t), t>0
\]
In analogy with the derivation of \eqref{3.8} this yields 
    
\begin{equation}\label{3.9}
    M_m=\gamma_2 \Psi(\gamma_2)
\end{equation}
with $i=-1$.

III 
\[
    u=P_m \Sigma_1^{j}\quad \textbf{on}\quad x=\Sigma_1(t), t>0,
\]
On insertion of the symmetry representation \eqref{2.2} there results $j=-1$ together with
\begin{equation}\label{3.10}
P_m=\gamma_1\Phi(\gamma_1)=L_m
\end{equation}

IV 
\[
    v=R_m \Sigma_2^{l}\quad \textbf{on}\quad x=\Sigma_2(t), t>0,
\]
Insertion of the symmetry representation \eqref{2.3} yields
\begin{equation}\label{3.11}
R_m=\gamma_2\Psi(\gamma_2)=M_m
\end{equation}
and $l=-1$.

In the preceding, it is required to determine $\Phi(\xi)$, $\Psi(\xi)$ via the coupled Ermakov-Ray-Reid/ Painlevé II system \eqref{2.8}-\eqref{2.9}. In \cite{rogers2016b} such were derived via the requirement that this system with $\alpha_I=\alpha_{II}=0$ has underlying Hamiltonian structure. In that case $S(\Psi/\Phi)$ and $T(\Phi/\Psi)$ were shown to admit the parametrisations
\begin{equation}\label{3.12}
S(\Psi/\Phi)=2(\Psi/\Phi)J(\Psi/\Phi)+(\Psi/\Phi)^{2}J'(\Psi/\Phi), 
\end{equation}

\[T(\Phi/\Psi)=-(\Psi/\Phi)^{2}J'(\Psi/\Phi)\]

It was established in \cite{rogers2016b} that the system their admits a Ermakov type invariant
\begin{equation}\label{3.13}
    I=\frac{1}{2}\left(\Phi\Psi_\xi-\Phi_\xi\Psi\right)^{2}+\frac{(\Phi^{2}+\Psi)^{2}}{\Phi^{2}}J(\Psi/\Phi)
\end{equation}
which was applied in the derivation of parametrisation of a class of exact solutions to the Ermakov-Ray-Reid/ Painlevé II system.
If the latter has underlying Hamiltonian structure, the Backlund transformation admitted by Painlevé II was applied iterativelyin \cite{rogers2016b} in this connection.

In conclusion, it is remarked that the specialisations
\[S=-1/3 (v/u)\quad, \quad\quad T=S=-1/3 (v/u)\] in \eqref{2.1} result in a novel two-component mKdV system with temporal modulation. In this case, symmetry reduction via \eqref{2.2} and \eqref{2.3} results in the hybrid Ermakov-Painlevé II system of a type  originally derived in \cite{rogers2014} by means of a wave packet representation admitted by a multi-component nonlinear Schrödinger system which incorporates de Broglie-Bohm type potential terms.

\subsection*{Acknowledgements}

The present work has been partially sponsored by the Project PIP No 11220220100532CO CONICET-UA, Rosario, Argentina and the Project O06-26CI2100 Universidad Austral Rosario Argentina.


\begin{thebibliography}{100} 


\bibitem{amster2015} Amster P and Rogers C, On a Ermakov-Painlevé II reduction in three-ion electrodiffusion. A Dirichlet boundary value problem, {\it Discrete and Continuous Dynamical Systems}, 2015, V. 35, 3277--3292.

\bibitem{calogero1985} Calogero F and Degasperis A, A symmetry approach to exactly solvable evolution equations, {\it J. Math. Phys.} 1985,  21, 1318–1325.


\bibitem{camassa1993} Camassa R and Holm D. An integrable shallow water equation with peaked solitons, {\it Phys. Rev. Lett.} 1993, V. 71, 1661--1664.

\bibitem{ermakov1880} Ermakov V P, Second-order differential equations: Conditions of complete integrability {\it Math. Proc. Camb. Phi. Soc.} Univ. Izy. Kiev 1880, 20, 1--25

\bibitem{flavin2008} Flavin J N and Rogers C, Upper estimates for a moving boundary problem for resonant nonlinear Schrödinger equations, {\it Studies in Applied Mathematics,} 2008,  121, 189-198.

\bibitem{fokas1980} Fokas A S, A symmetry approach to exactly solvable evolution equations, {\it J. Math. Phys.} 1980, 21, 1318–1325.



\bibitem{fokas1998} Fokas A S and Tanveer S, A Hele-Shaw problem and the second Painlevé transcendent, {\it Math. Proc. Camb. Phi. Soc.}, 1998, V. 124, 169--191.


\bibitem{kingston1982}  Kingston J G and Rogers C, Reciprocal Bäcklund transformations of conservation laws, {\it Phys. Lett.}, 1982, V. 92A, 261--264.

\bibitem{lee2007} Lee J H, Pashaev O K, Rogers C and Schief W K, The resonant nonlinear Schrödinger equation in cold plasma physics. application of Bäcklund-Darboux transformations and superposition principles, {\it Journal of Plasma Physics}, 2007, 73, 257-272.

\bibitem{lukashovich1971} Lukashevich N A, The second Painlevé equation, {\it Differential Equations},1971,  7, 853-854.
\bibitem{ray1980} Ray J R, Nonlinear superposition law for generalized Ermakov systems, {\it Physics Letters A}, 1980, 78, 4-6 .

\bibitem{reid1980} Reid J L and Ray J R, Ermakov systems, nonlinear superposition, and solutions of nonlinear equations of motion, {\it J. Math. Phys.}, 1980, 21, 1583--1587.

\bibitem{rogers1989a} Rogers C and Ramgulam U, A nonlinear superposition principle and Lie group invariance. Application in rotating shallow water theory, {\it Int. J. Nonlinear Mech.}, 1989, V.24, 229--235.

\bibitem{rogers1989b} Rogers C and Ames W F, Nonlinear boundary value problems in science and engineering,  Academic Press, New York, 1989.

\bibitem{rogers1989c} Rogers C, Elliptic warm core theory: the pulsrodon, {\it Phys. Lett.}, 1989, A 138, 267--273.

\bibitem{rogers1993} Rogers C, Hoenselaers C and Ray J R, On (2+1)-dimensional Ermakov systems, {\it J. Phys. A: Math. Gen.} 1993, 26,  2625--2633

\bibitem{rogers1996} Rogers C,  Schief W K ,
Multi-component Ermakov systems: structure and linearization, {\it
Journal of Mathematical Analysis and Applications}, 1996, 198, 194-220.

\bibitem{rogers1999} Rogers C,  Schief W K , The resonant nonlinear Schrödinger equation via an integrable capillarity model,
  {\it 
Journal of Mathematical Analysis and Applications}, 1999, 198, 194-220.

\bibitem{rogers2009} Rogers C,  Yip L P and Chow K W, A resonant Davey-Stewartson capillarity model system: solitonic generation, {\it International Journal of Nonlinear Sciences and Numerical Simulation}, 2009, 10, 397--405.

\bibitem{rogers2010a} Rogers C,  Malomed B,  Chow K and An H, Ermakov–Ray–Reid systems in nonlinear optics, {\it . Phys. A: Math. Theor. } 2010,  43,  455214 (15 pp).

\bibitem{rogers2010b} Rogers C,  A Ermakov–Ray–Reid reduction in 2+1-dimensional magnetogasdynamics, {\it .Proceedings, 5th International Workshop in Group Analysis of Differential Equations and Integrable Systems:  Cyprus}, 2010, 164-177

\bibitem{rogers2010c} Rogers C and An H, Ermakov–Ray–Reid systems in (2+1)-dimensional rotating shallow water theory, {\it Studies in Applied Mathematics}, 2010, 125: 275-299.

\bibitem{rogers2011} Rogers C and An H, On a (2+1)-dimensional Madelung system with logarithmic and de Broglie- Bohm quantum potentials: Ermakov reduction, {\it Phys. Scr.} (2011), 84 045004 (7 pp)

\bibitem{rogers2012} Rogers C, Malomed B and An H, Ermakov-Ray-Reid reductions of variational approximations in nonlinear optics, {\it Stud. Appl. Math}, 2012, V. 129, 389--413.

\bibitem{rogers2013} Rogers C, On a Coupled Nonlinear Schrödinger System: A Ermakov Connection, {\it Stud. Appl. Math.}, 2013, V. 132, 247--265.

\bibitem{rogers2014} Rogers C, A novel Ermakov-Painlevé II System. N+1-dimensional coupled NLS and elastodynamic reductions, {\it Stud. Appl. Math.}, 2014, V. 133, 214--231.


\bibitem{rogers2015} Rogers C, Moving Boundary Problems for the Harry Dym equation and its reciprocal associates, {\it Zeit angew Math. Phys.}, 2015, V.66, 3205--3220.

\bibitem{rogers2016} Rogers C, On a class of moving boundary problems for the potential mkdV equation: conjugation of Bäcklund and reciprocal transformations, {\it Ricerche di Matematica}, 2016, 65, 563–577



\bibitem{rogers2016b} Rogers C and  Schief W K, On Ermakov-Painlevé II systems. Integrable reduction, {\it Meccanica}, 2016, V.51, 2967--2974.



\bibitem {rogers2017b} Rogers C, Moving boundary problems for an extended Dym equation. Reciprocal connections, {\it Meccanica}, 2017, V.52, 3531--3540.


\bibitem{rogers2017c} Rogers C and  Clarkson P A, Ermakov-Painlevé II symmetry reduction of a Korteweg capillarity system, {\it Symmetry, Integrability and Geometry: Methods and Applications}, 2017, 13, 018.

\bibitem{rogers2018a} Rogers C and  Clarkson P A, Ermakov-Painlevé II reduction in cold plasma physics. Application of a Bäcklund transformation, {\it J. Nonlinear Mathematical Physics}, 2018, V.25, 247--261.

\bibitem{rogers2018b} Rogers C and Schief  W K, Ermakov-type systems in nonlinear physics and continuum mechanics, {\it Nonlinear Systems and Their Remarkable Mathematical Structures}, Ed. Norbert Euler, CRC Press, 2018, 541--576.

\bibitem{rogers2018} Rogers, C, Multi-component Ermakov and non-autonomous many-body system connections, {\it Ricerche di Matematica} 2018, 67,  803 - 815.

\bibitem{rogers2019b} Rogers C, On a canonical nine-body problem. Integrable Ermakov decomposition via reciprocal transformations {\it J Nonlinear Math. Phys.} 2019, 26, 98–106. 

\bibitem{rogers2022} Rogers C, Moving boundary problems for a canonical member of the WKI inverse scattering scheme: conjugation of a reciprocal and Möbius transformation, {\it Physica Scripta }, 2022, V.97, 095207.

\bibitem{rogers2023} Rogers C, On mKdV and associated classes of moving boundary problems: reciprocal connections, {\it Meccanica} 2023, V.58, 1633--1640.

\bibitem{rogers2025a} Rogers C, On Korteweg–de Vries and associated reciprocal moving boundary problems, {\it Z. angew. Math. Phys.}, 2025, V.76, 33. 

\bibitem{rogers2025b} Rogers C, On moving boundary problems for the solitonic Gardner equation. a reciprocally associated class, {\it Z. angew. Math. Phys.}, 2025, 76:186.

\bibitem{rogers2025c} Rogers C,  Briozzo A C, Moving boundary problems with Ermakov symmetry reduction: nonlinear superposition principle and reciprocal transformation applications, {\it Open Communications in Nonlinear Mathematical Physics}, 2025, 5, 81-91.


\bibitem{rogers2025d} Rogers C,  Briozzo A C, Moving boundary problems for a novel extended mKdV equation. Application of Ermakov-Painlevé II symmetry reduction, {\it Open Communications in Nonlinear Mathematical Physics}, 2025, 5, 116-130.


\bibitem{saffman1958} Saffman P G and  Taylor G I, The penetration of a fluid into a porous medium or Hele-Shaw cell containing a more viscous liquid, {\it Proceedings of the Royal Society of London}, 1958 A 245, 312--329.

\bibitem{shahinpoor1971} Shahinpoor M and Nowinski J L, Exact solutions to the problem of forced large amplitude radial oscillations in a thin hyperelastic tube, {\it Int. J. Nonlinear Mech.} 1971 V.6, 193-207.

\bibitem{schief1999} Schief W K and Rogers C, Binormal motion of curves of constant curvature and torsion. Generation of soliton surfaces, {\it Proc. Roy. Soc. London}, 1999 455 A, 3163--3188.  



\bibitem{vassiliou2001}  Vassiliou P J, Harry Dym equation in Encyclopaedia of Mathematics, Springer, 2001.




\end{thebibliography}
\end{document}